\documentclass[11pt]{article}

\usepackage{amsmath,graphicx}
\usepackage[T1]{fontenc} 
\usepackage[latin1]{inputenc}
\usepackage[frenchb,english]{babel}

\providecommand*{\fullref}[1]{\hyperref[{#1}]{\cref*{#1}. \nameref*{#1}}}
\providecommand*{\Fullref}[1]{\hyperref[{#1}]{\Cref*{#1}. \nameref*{#1}}}

\usepackage{amssymb} 
\usepackage{amsmath} 
\usepackage{amsfonts} 
\usepackage{mathrsfs}
\usepackage{amsthm}
\usepackage{kmath} 			
\usepackage[noadjust]{cite}       
\usepackage{soul}		         

\usepackage{graphicx}
\usepackage{color}
\usepackage{pgfplots}
\usepackage{epsfig}
\usepackage{epstopdf}

\def\0{{\mathbf 0}}
\def\1{{\mathbf 1}}

\def\G{{\mathbf G}}

\def\I{{\mathbf I}}

\def\U{{\mathbf U}}

\def\X{{\mathbf X}}

\def\aa{{\boldsymbol{a}}}

\def\dd{{\boldsymbol d}}

\def\hh{{\boldsymbol h}}

\def\vv{{\boldsymbol v}}

\def\zz{{\boldsymbol z}}

\def\Dc{{\mathcal D}}

\def\Hc{{\mathcal H}}

\def\Sc{{\mathcal S}}

\def\Rbb{{\mathbb R}}

\def\ie{\textit{i.e.},\xspace}
\def\eg{\textit{e.g.},\xspace}

\newkmacro\sign[1][{\cdot}]{\mathrm{sign}\left({#1}\right)} 
\newkmacro\interior[1][{\cdot}]{\mathrm{int}\left({#1}\right)}
\newkmacro\bd[1][{\cdot}]{\mathrm{bd}\left({#1}\right)}
\newkmacro\conv[1][{\cdot}]{\mathrm{conv}\left({#1}\right)}
\newkmacro\card[1][{\cdot}]{\mathrm{card}\left({#1}\right)}
\newkmacro\spa[1][{\cdot}]{\mathrm{span}\left({#1}\right)}
\newkmacro\nul[1][{\cdot}]{\mathrm{null}\left({#1}\right)}
\newkmacro\deter[1][{\cdot}]{\mathrm{det}\left({#1}\right)}
\newkmacro\tr[1][{\cdot}]{\mathrm{tr}\left({#1}\right)}

\newtheorem{theorem}{Theorem}

\newtheorem{example}{Example}

\def\Sh{{\Hc}} 				
\def\vech{{\hh}} 		         
\def\gram{\G}

\def\Spwidth{{\epsilon}}		

\def\dimSps{{n}} 		        	
\def\Sps{{V}}				
\def\Spswidth{{\hat{\Spwidth}}}      
\def\Spsbasis{{\vv}} 		        
\def\Spsp{{\Sps}^\perp}		

\def\idxss{{k}}

\def\dimSo{{m}} 		        
\def\So{{W}} 				

\def\SoPG{{Z}}
\def\SoPGbasis{{\zz}}
\def\rieszr{{\aa}} 				

\newkmacro{\dist}[2][]{\mathrm{dist}(#1,#2)}
\newkmacro{\scap}[2][]{\kangle{#1,#2}}
\newkfunc{\projector}{P}

\def\svdecr{{\sigma}} 		                

\def\aop{{a}}

\def\lop{{b}}

\def\zvar{{\vech}}

\def\ztrue{\zvar^\star}
\def\zrom{{\hat{\zvar}}}
\def\zromPG{{\zrom_{\mathrm{PG}}}}

\def\zrommult{{\zrom_{\mathrm{MS}}}}
\def\tmpvar{{\eta}}

\def\feasSd{{\Dc}}
\def\feasd{{\dd}}
\def\slice{\Sc}

\title{A Mathematical Characterization of the Performance of the ``Multi-Slice'' Projector}
\author{C. Herzet, M. Diallo, P. Héas \\
\small{INRIA Rennes, Campus de Beaulieu, 35000 Rennes, France}}
\date{}

\begin{document}
%
\maketitle

\abstract
We consider an enhanced version of the well-kwown ``Petrov-Galerkin'' projection in Hilbert spaces. The proposed procedure, dubbed ``multi-slice'' projector, exploits the fact that the sought solution belongs to the intersection of several high-dimensional slices. This setup is for example of interest in model-order reduction where this type of prior may be computed off-line.  In this note, we provide a mathematical characterization of the performance achievable by the multi-slice projector and compare the latter with the results holding in the Petrov-Galerkin setup. In particular, we illustrate the superiority of the multi-slice approach in certain situations. \\

Nous considérons une version améliorée de la projection de ``Petrov-Galerkin'' dans un espace de Hilbert. La procédure proposée, appelée ``projecteur multi-tranches'', exploite le fait que la solution recherchée appartient à l'intersection de plusieurs tranches de hautes dimensions. Dans cette note, nous fournissons une caractérisation mathématique des performances atteignables par le projecteur ``multi-tranches'' et comparons les résultats obtenus à ceux existants dans le contexte des projections de Petrov-Galerkin. Nous illustrons ainsi la supériorité de l'approche multi-tranches dans certaines situations.


\section{Introduction}

Let $\Sh$ be a Hilbert space with inner product $\scap[\cdot]{\cdot}$ and induced norm $\kvvbar{\cdot}$. We focus on the following variational formulation:\vspace{0.2cm}
\begin{align}
\mbox{Find $\ztrue\in\Sh$ such that } \ \aop(\ztrue, \vech) = \lop(\vech)\quad \forall \vech\in\Sh, \label{eq:probleminit}
\end{align}
where $\aop: \Sh\times\Sh \rightarrow \Rbb$  is a bilinear operator and $\lop: \Sh\rightarrow \Rbb$ a linear operator.~Problem \eqref{eq:probleminit} is quite common (it appears for example in the weak formulation of elliptic partial differential equations) and has therefore been well-studied in the literature.~In particular, it has a unique solution under mild conditions, see Lax-Milgram's and Necas Theorems in \cite[Theorems 2.1 and 2.2]{Quarteroni2016Reduced}. \vspace{0.2cm}

Unfortunately, solving \eqref{eq:probleminit} is generaly an intractable problem. A popular alternative to compute an approximation of \eqref{eq:probleminit} is known as ``Petrov-Galerkin'' projection.~Formally, this approach consists of approximating \eqref{eq:probleminit} by the following problem:\vspace{0.2cm}
\begin{align}
\mbox{Find $\zromPG\in\Sps_\dimSps$ such that } \ \aop(\zromPG, \vech) = \lop(\vech)\quad \forall \vech\in\SoPG_\dimSo \label{eq:problemPG}\\[-0.4cm]\nonumber
\end{align}
where $\Sps_\dimSps\subset\Sh$ is a linear subspace of dimension $\dimSps$ and $\SoPG_\dimSo\subset\Sh$  is a linear subspace of dimension $\dimSo\geq \dimSps$. Since the dimension of $\Sps_\dimSps$ and $\SoPG_\dimSo$ are finite, \eqref{eq:problemPG} admits a simple algebraic solution under mild conditions.  In the literature of model reduction (see \eg \cite{Quarteroni2016Reduced}), Petrov-Galerkin approximation is at the core of the family of ``projection-based'' reduced models. \vspace{0.2cm}

In this note we elaborate on an alternative projection procedure exploiting several approximation subspaces.~Indeed, in the context of model-order reduction, standard strategies to evaluate a good approximation subspace $\Sps_\dimSps$, \eg reduced basis \cite{Quarteroni2016Reduced} or proper orthogonal decomposition \cite{Berkooz1993Proper}, typically generate a sequence of subspaces $\kbrace{\Sps_\idxss}_{\idxss=0}^\dimSps$ and positive scalars $\kbrace{\Spswidth_\idxss}_{\idxss=0}^\dimSps$ such that\vspace{0.2cm}
\begin{align}
\Sps_0\subset \Sps_1\subset \ldots \subset \Sps_\dimSps
\end{align}
and
\begin{align}
\dist[\ztrue]{\Sps_\idxss}\leq \Spswidth_\idxss,\quad \idxss = 0\ldots\dimSps. \label{eq:msingr1}
\end{align}
Clearly, \eqref{eq:msingr1} provides some useful information about the location of $\ztrue$ in $\Sh$ since it restrains the latter to belong to the intersection of a set of low dimensional slices, \ie\vspace{0.1cm}
\begin{align}
\ztrue \in \cap_{\idxss=0}^\dimSps \slice_{\idxss}, \label{eq:priorslices}
\end{align}
where
\begin{align}
\slice_{\idxss} &= \kbrace{ \vech : \dist[\vech]{\Sps_\idxss}\leq \Spswidth_\idxss}, \quad \idxss = 0\ldots\dimSps.
\end{align}
 In standard Petrov-Galerkin projection \eqref{eq:problemPG}, only $\Sps_\dimSps$ is used and the additional information provided by \eqref{eq:priorslices} is discarded. In this work, we consider a simple methodology to exploit the latter additional information into the projection process.
More specifically, we focus on the following optimization problem\footnote{In this note we assume that constraints are available $\forall \idxss\in\kbrace{1\ldots\dimSps}$. All the derivations presented in this paper may nevertheless be easily extended to the case where constraints in \eqref{eq:MSdec} are only available for \textit{some} $\idxss\in\kbrace{1\ldots\dimSps}$. } 
\begin{align}
\mbox{Find }\, \zrommult
\in & \kargmin_{\zvar\in\Sps_\dimSps} \sum_{j=1}^\dimSo \kparen{\lop(\SoPGbasis_j)- \aop(\zvar,\SoPGbasis_j)}^2 \label{eq:MSdec}\\
&\mbox{subject to } \dist[\vech]{\Sps_\idxss}\leq \Spswidth_\idxss,\quad \idxss = 0\ldots\dimSps,\nonumber
\end{align}
which can be seen as an extension of the standard Petrov-Galerkin approach.~In particular, the constraints in \eqref{eq:MSdec} exploit the prior information \eqref{eq:msingr1} into the projection process: each constraint imposes that the solution belongs to some $\idxss$-dimensional slice $\slice_{\idxss}$.~Hence, in the sequel, we will dub this procedure as ``multi-slice'' projection. \vspace{0.2cm}

The practical interest of the multi-slice approach has already been emphasized in several contributions.~In \cite{Herzet:2017fb, Herzet:2018qw} we presented some applications of the multi-slice decoder to the problem of model-order reduction of parametric partial differential equations.~In \cite{JP-Argaud:2017qm} and \cite{2017arXiv171003569F}, the authors showed that multi-slice decoder can be of interest to enhance the performance of  the ``empirical interpolation method'' or the simulation of Navier-Stokes equations.~``Multi-slice'' prior information of the form \eqref{eq:priorslices} has also been considered in \cite{Binev2015Data} for data assimilation. However, in the latter contribution, the decoder considered by the authors differs from \eqref{eq:MSdec} since the solution is no longer constrained to belong to the low-dimensional subspace $\Sps_\dimSps$.  \vspace{0.2cm}

In this note we provide a mathematical characterization of the performance achievable by the multi-slice decoder \eqref{eq:MSdec}. More specifically, we derive an ``instance optimality property'' relating the projection error $\|\zrommult-\ztrue\|$ to the distance between $\ztrue$ and the different approximation subspaces $\Sps_\idxss$. Our result is presented in  Theorem~\ref{eq:UBMS} in the next section.


\section{Performance guarantees}\label{sec:IOP}

One of the reasons which has ensured the success of Petrov-Galerkin projection is the existence of strong theoretical guarantees, \eg Cea's Lemma~\cite[Lemma 2.2]{Quarteroni2016Reduced} or the  Babuska's Theorem \cite[Theorem 2.3]{Quarteroni2016Reduced}. In this section we derive a similar result for the multi-slice decoder \eqref{eq:MSdec}. The standard result associated to Petrov-Galerkin projection is recalled in Theorem~\ref{th:UBPG} whereas our characterization of the multi-slice decoder \eqref{eq:MSdec} is presented in Theorem~\ref{eq:UBMS}. 
We conclude this section by providing two examples in which the multi-slice projector leads to better guarantees of reconstruction than the standard Petrov-Galerkin approach. \vspace{0.2cm}

We first introduce some quantities of interest.~First, we let $\kbrace{\Spsbasis_j}_{j=1}^\dimSps$ and $\kbrace{\SoPGbasis_j}_{j=1}^\dimSo$ be orthonormal bases (ONBs) of  the subspaces $\Sps_\dimSps$ and $\SoPG_\dimSo$, respectively.~We define $\kbrace{\rieszr_j}_{j=1}^\dimSo$ as the Riesz's representers of $\kbrace{\aop(\cdot,\SoPGbasis_j)}_{j=1}^\dimSo$.  We denote by $\kbrace{\svdecr_j}_{j=1}^\dimSps$ the set of singular values (sorted in their decreasing order of magnitude) of the Gram matrix
\begin{align}
\gram = [\scap[\rieszr_i]{\Spsbasis_j}]_{i,j}\in\Rbb^{\dimSo\times\dimSps}. \label{eq:Gram}
\end{align}
With these notations, the well-known Babuska's theorem (in a Hilbert space) can be formulated as follows:\vspace{0.2cm}

\begin{theorem}[Babuska's Theorem]\label{th:UBPG} If $\svdecr_\dimSps>0$ then the solution of \eqref{eq:problemPG} is unique and satisfies
\begin{align}
\kvvbar{\ztrue-\zromPG} \leq \frac{\svdecr_1}{\svdecr_\dimSps} \dist[\ztrue]{\Sps_\dimSps}.  \label{eq:IOPSS}
\end{align}
\end{theorem}
\vspace{0.2cm}
See for example \cite{Xu:2003lr} for a proof of this result. 
Hereafter we provide a similar characterization of the performance of the multi-slice projector~\eqref{eq:MSdec}. In order to state our result we need to introduce the following quantities. We first define the short-hand notations\footnote{$\Spwidth_\idxss$ thus represents the true distance from $\ztrue$ to $\Sps_\idxss$. We note that this quantity is usually unknown to the practitioner.  This is in contrast which $\Spswidth_\idxss$ which represents the prior information available to the practitioner but is only an upper bound on $\Spwidth_\idxss$.}
\begin{align}
\Spwidth_\idxss = \dist[\ztrue]{\Sps_\idxss}, 
\end{align}
and 
\begin{align}
\gamma = \sup_{\vech\in\Spsp_\dimSps, \kvvbar{\vech}=1} \kparen{\sum_{j=1}^\dimSo \scap[\rieszr_j]{\vech}^2}^{\frac{1}{2}}.\label{eq:defgamma}
\end{align}
%
Moreover, we define
\begin{align}
\delta_j = \sum_{\idxss=1}^\dimSps \kvbar{x_{\idxss j}} (\Spswidth_{\idxss-1}+\Spwidth_{\idxss-1}),
\end{align}
where $x_{\idxss j}$ are the elements of the matrix $\X$ appearing in the singular value decomposition of $\G$, that is $\G = \U \Lambda \ktranspose{\X}$, where $\U\in\Rbb^{\dimSo\times\dimSo}$, $\X\in\Rbb^{\dimSps\times\dimSps}$ are orthogonal matrices and  $\Lambda\in\Rbb^{\dimSo\times\dimSps}$ is the diagonal matrix of singular values $\kbrace{\svdecr_j}_{j=1}^\dimSps$. 

\vspace{0.2cm}
Using these notations, our result reads:
\begin{theorem}\label{eq:UBMS}
Let $\ztrue$ be a solution of \eqref{eq:probleminit} verifying  \eqref{eq:priorslices}. 
Then any solution $\zrommult$ of \eqref{eq:MSdec} verifies 
\begin{align}\label{eq:WCUB}
\kvvbar{\ztrue - \zrommult} \leq 
\left\{
\begin{array}{ll}
\kparen{\sum_{j=\ell+1}^\dimSps \delta_j^2+ \rho\, \delta_\ell^2+ \Spwidth_\dimSps^2}^{\frac{1}{2}}& \mbox{if $\sum_{j=1}^{\dimSps} \svdecr_j^2 \delta_j^2\geq 4 \gamma^2 \Spwidth_\dimSps^2$},\\
\kparen{\sum_{j=1}^{\dimSps} \delta_j^2+ \Spwidth_\dimSps^2}^{\frac{1}{2}}& \mbox{otherwise,}
\end{array}
\right.
\end{align}
where $\ell$ is the largest integer such that
\begin{align}
\sum_{j=\ell}^{\dimSps} \svdecr_j^2 \delta_j^2\geq 4 \gamma^2 \Spwidth_\dimSps^2,
\end{align}
and $\rho\in[0,1]$ is defined as
\begin{align}
\rho \svdecr_\ell^2 \delta_\ell^2 +\sum_{j=\ell+1}^{\dimSps} \svdecr_j^2 \delta_j^2= 4 \gamma^2 \Spwidth_\dimSps^2. \label{eq:WCUB2}
\end{align}
Moreover, if $\svdecr_\dimSps>0$, \eqref{eq:MSdec} admits a unique solution. 
\end{theorem}

A proof of Theorem \ref{eq:UBMS} is detailed in Section~\ref{sec:proofMainTh}. \vspace{0.2cm}


We conclude this section by particularizing the results stated in Theorems~\ref{th:UBPG} and \ref{eq:UBMS} to different setups. In particular, we emphasize two situations\footnote{The two setups considered below correspond to those exposed in \cite[Section 3.2]{Binev2015Data}.} where the multi-slice projection has much better reconstruction guarantees than its Petrov-Galerkin counterpart. In order to ease the comparison between the bounds stated in Theorems \ref{th:UBPG} and \ref{eq:UBMS}, we consider the case where $\kbrace{\rieszr_j}_{j=1}^\dimSo$ is an ONB. We note that in such a case, we have $\svdecr_1\leq 1$ and $\gamma\leq 1$. \\ 

\begin{example} 
We first assume that $\X=\I_\dimSps$ in the singular-value decomposition of $\G$. 
We set  $\Spswidth_j = \Spwidth_j$ and assume that
\begin{align}
\Spwidth_j &= \left\{
\begin{array}{ll}
1 & j = 0 \ldots \dimSps-3,\\
\Spwidth^{\frac{1}{2}} & j = \dimSps-2, \dimSps-1,\\
\Spwidth & j=\dimSps, 
\end{array}
\right. 
\end{align}
for some $\Spwidth\ll 1$. 
Moreover, we let
\begin{align}
\svdecr_j &= \left\{
\begin{array}{ll}
1 & j = 1 \ldots \dimSps-3,\\
\Spwidth^{\frac{1}{2}} & j = \dimSps-2, \dimSps-1,\\
\Spwidth & j=\dimSps. 
\end{array}
\right. 
\end{align}
In this setup, the upper bound \eqref{eq:IOPSS} of Theorem \ref{th:UBPG} becomes:
\begin{align}
\kvvbar{\zromPG -  \ztrue}\leq \svdecr_\dimSps^{-1}\ \dist[\ztrue]{\Sps_\dimSps}
&= \Spwidth^{-1} \Spwidth = 1. \label{eq:particularBTh1}
\end{align}
On the other hand, because $\X = \I$, we have  
\begin{align}
\delta_j = \Spswidth_{j-1}+\Spwidth_{j-1} = 2 \Spwidth_{j-1}.
\end{align}
The index $\ell$ appearing in Theorem \ref{eq:UBMS} is smaller or equal to $\dimSps-1$ since
\begin{align}
\svdecr^2_\dimSps \delta_\dimSps^2 &= \svdecr^2_\dimSps (2 \Spwidth_{\dimSps-1})^2 = 4 \Spwidth^3 \ll 4 \Spwidth^2, \nonumber\\
\svdecr^2_{\dimSps-1} \delta_{\dimSps-1}^2 &= \svdecr^2_{\dimSps-1} (2 \Spwidth_{\dimSps-2})^2 = 4 \Spwidth^2, \nonumber
\end{align}
and thus
\begin{align}
\svdecr^2_{\dimSps-1} \delta_{\dimSps-1}^2+ \svdecr^2_\dimSps \delta_\dimSps^2\geq 4 \Spwidth^2 \geq 4 \gamma^2 \Spwidth^2
\end{align}
since $\gamma \leq 1$. The upper bound in Theorem \ref{eq:UBMS} becomes
\begin{align}
\kvvbar{\ztrue - \zrommult} 
&\leq  \kparen{\delta_{\dimSps-1}^2+ \delta_\dimSps^2 + \Spwidth_\dimSps^2}^{\frac{1}{2}},\nonumber\\
&= \kparen{4 \Spwidth+ 4 \Spwidth + \Spwidth^2}^{\frac{1}{2}},\nonumber\\
&\leq 3 \Spwidth^{\frac{1}{2}}. \label{eq:particularBTh2a}
\end{align}
Hence the bound in the multi-slice setup \eqref{eq:particularBTh2a} can be arbitrarily small as compared to \eqref{eq:particularBTh1} when $\Spwidth \rightarrow 0$. \\
\end{example}

\begin{example}
We now consider $\X= \dimSps^{-\frac{1}{2}}\, \1_{\dimSps\times\dimSps}$ where $\1_{\dimSps\times\dimSps}$ is an $\dimSps\times\dimSps$ matrix of $1$'s.   
%
We set  $\Spswidth_j = \Spwidth_j$ and assume that
\begin{align}
\Spwidth_j &= \left\{
\begin{array}{ll}
\frac{1}{2} & j = 0,\\
\frac{1}{2(\dimSps-1)} & j = 1 \ldots \dimSps-1,\\
\Spwidth & j=\dimSps, 
\end{array}
\right. 
\end{align}
for some $\Spwidth\ll \dimSps^{-1}$ (Note that we must have: $\Spwidth \leq  \frac{1}{2(\dimSps-1)}$ by definition). 
Moreover, we let
\begin{align}\label{eq:defsvex2}
\svdecr_j &= \left\{
\begin{array}{ll}
\svdecr & j = 1 \ldots \dimSps-1,\\
\Spwidth^2 & j=\dimSps, 
\end{array}
\right. 
\end{align}
for some $1\geq \svdecr>\Spwidth$ whose value will be specified below. 

With these choices, the upper bound \eqref{eq:IOPSS} of Theorem \ref{th:UBPG} becomes:
\begin{align}
\kvvbar{\zromPG -  \ztrue}\leq \svdecr_\dimSps^{-1}\ \dist[\ztrue]{\Sps_\dimSps}
&= \Spwidth^{-2} \Spwidth = \Spwidth^{-1}. \label{eq:particularBTh1_ex2}
\end{align}
On the other hand, we have 
\begin{align}
\delta_j 
&= \sum_{\idxss=1}^\dimSps \kvbar{x_{\idxss j}} (\Spswidth_{\idxss-1}+\Spwidth_{\idxss-1}),\nonumber\\
&= 2 \dimSps^{-\frac{1}{2}}  \sum_{\idxss=1}^\dimSps \Spwidth_{\idxss-1},\nonumber\\
&= 2 \dimSps^{-\frac{1}{2}}. 
\end{align}
By choosing $\svdecr$ such that (we remind the reader that $\svdecr_{\dimSps-1}=\svdecr$ by definition \eqref{eq:defsvex2})
\begin{align}
 \svdecr_{\dimSps-1}^2\delta_{\dimSps-1}^2+\svdecr_{\dimSps}^2\delta_{\dimSps}^2 = 4 \Spwidth^2, 
\end{align}
we obtain that index $\ell$ appearing in Theorem \ref{eq:UBMS} is smaller or equal to $\dimSps-1$ since $\gamma\leq 1$. 
The upper bound in Theorem \ref{eq:UBMS} then reads
\begin{align}
\kvvbar{\ztrue - \zrommult} 
&\leq  \kparen{\delta_{\dimSps-1}^2+ \delta_\dimSps^2 + \Spwidth_\dimSps^2}^{\frac{1}{2}},\nonumber\\
&= \kparen{4 \dimSps^{-1}+ 4 \dimSps^{-1} + \Spwidth^2}^{\frac{1}{2}},\nonumber\\
&\leq 3 \dimSps^{-\frac{1}{2}}, \label{eq:particularBTh2}
\end{align}
where the last inequality follows from our initial assumption $\Spwidth\ll \dimSps^{-1}$. \\
\end{example}

\section{Proof of Theorem \ref{eq:UBMS}}\label{sec:proofMainTh}

In this section, we provide a proof of the result stated in Theorem \ref{eq:UBMS}. 
 We first note that problem~\eqref{eq:MSdec} is equivalent to finding the minimum of a quadratic function over a closed bounded subset of $\Sps_\dimSps$. A minimizer thus always exists. Moreover, the unicity of the minimizer stated at the end of Theorem~\ref{eq:UBMS} follows from the strict convexity of the cost function when $\svdecr_\dimSps>0$. \vspace{0.2cm}

In the rest of this section, we thus mainly focus on the derivation of the upper bound \eqref{eq:WCUB}.~Our proof is based on the following steps. First, since $\zrommult\in\Sps_\dimSps$, we have that 
\begin{align}
\kvvbar{\ztrue - \zrommult}^2 
&= \kvvbar{\projector[\Sps_\dimSps](\ztrue) - \zrommult}^2 + \kvvbar{\projector[\Sps_\dimSps][\perp](\ztrue)}^2, \nonumber\\
&= \kvvbar{\projector[\Sps_\dimSps](\ztrue) - \zrommult}^2 + \Spwidth_\dimSps^2, \label{eq:firstubound}
\end{align}
where $\projector[\Sps_\dimSps](\cdot)$ (resp. $\projector[\Sps_\dimSps][\perp](\cdot)$) denotes the orthogonal projector onto $\Sps_\dimSps$ (resp. $\Sps_\dimSps^\perp$). 
We then derive an upper bound on $\|\projector[\Sps_\dimSps](\ztrue) - \zrommult\|^2$ as follows:
\begin{itemize}
\item We identify a set $\feasSd$ such that $\projector[\Sps_\dimSps](\ztrue)-\zrommult\in\feasSd$ in Section~\ref{sec:defD}. We then have $\|\projector[\Sps_\dimSps](\ztrue)-\zrommult\|^2\leq \sup_{\feasd\in\feasSd}\kvvbar{\feasd}^2$.
\item We derive the analytical expression of $\sup_{\feasd\in\feasSd}\kvvbar{\feasd}^2$ as a function of the parameters $\kbrace{\Spwidth_\idxss}_{\idxss=1}^\dimSps$, $\kbrace{\Spswidth_\idxss}_{\idxss=1}^\dimSps$ and $\kbrace{\svdecr_\idxss}_{\idxss=1}^\dimSps$.
\end{itemize}
Combining these results, we obtain \eqref{eq:WCUB}-\eqref{eq:WCUB2}.

\subsection{Definition of $\feasSd$}\label{sec:defD}

We express $\feasSd$ as the intersection of two sets $\feasSd_1$ and $\feasSd_2$ that we define in Sections \ref{sec:defD1} and \ref{sec:defD2} respectively. In order to properly define these quantities, we introduce some particular  ONBs for $\Sps_\dimSps$ and $\So_{\dimSo}=\spa[\kbrace{\rieszr_j}_{j=1}^\dimSo]$ in Section~\ref{sec:BasesFavorables}.  

\subsubsection{Some particular bases for $\Sps_\dimSps$ and $\So_{\dimSo}$}\label{sec:BasesFavorables}

%

Let 
\begin{align}
\gram = \U \Lambda \ktranspose{\X} \label{eq:svd_Gr}
\end{align}
be the singular value decomposition of the Gram matrix defined in \eqref{eq:Gram}, 
 where $\U\in\Rbb^{\dimSo\times\dimSo}$ and $\X\in\Rbb^{\dimSps\times\dimSps}$ are orthonormal matrices and $\Lambda\in\Rbb^{\dimSo\times\dimSps}$ is the diagonal matrix of singular values. We  denote by $\kbrace{\svdecr_j}_{j=1}^\dimSps$ the set of singular values of $\G$ sorted in their decreasing order of magnitude.  \vspace{0.2cm}

We define the following bases for $\Sps_\dimSps$ and $\So_{\dimSo}$:
\begin{align}
\Spsbasis_j^{*} &= \sum_{i=1}^\dimSps x_{ij} \Spsbasis_i,\\
\rieszr_j^{*} &= \sum_{i=1}^\dimSo u_{ij} \rieszr_i,
\end{align}
where $\U\in\Rbb^{\dimSo\times\dimSo}$ and $\X\in\Rbb^{\dimSps\times\dimSps}$ are the orthonormal matrices appearing in \eqref{eq:svd_Gr}. We note that $\kbrace{\Spsbasis_j^*}_{j=1}^\dimSps$ is an ONB whereas $\kbrace{\rieszr_j^*}_{j=1}^\dimSo$ is not necessarily orthonormal. By definition, $\kbrace{\Spsbasis_j^*}_{j=1}^\dimSps$ and $\kbrace{\rieszr_j^*}_{j=1}^\dimSo$ enjoy the following desirable property:
\begin{align}
\scap[\rieszr_i^*]{\Spsbasis_j^*} &= \left\{
\begin{array}{ll}
\svdecr_j & \mbox{if $i=j$} \\
0 & \mbox{otherwise.}
\end{array}
\right. \label{eq:orthproperty}
\end{align} 
\vspace{0.2cm}

\subsubsection{Definition of $\feasSd_1$}\label{sec:defD1}

Let us define $\feasSd_1$ as 
\begin{align}
\feasSd_1 = \kbrace{\feasd = \sum_{j=1}^{\dimSps} \beta_j \Spsbasis_j^* : \sum_{j=1}^{\dimSps} \svdecr_j^2 \beta_j^2\leq 4 \gamma^2 \Spwidth_\dimSps^2},
\end{align}
where $\gamma$ is defined in \eqref{eq:defgamma}.
 We show hereafter that $\projector[\Sps_\dimSps](\ztrue)-\zrommult\in\feasSd_1$.\vspace{0.2cm}

Let us first consider the intermediate set 
\begin{align}
\Sc&= \kbrace{\vech : f(\vech) \leq \gamma^2 \Spwidth_\dimSps^2}, 
\end{align}
where $f(\vech) = \sum_{j=1}^\dimSo \kparen{\lop(\SoPGbasis_j)- \aop(\zvar,\SoPGbasis_j)}^2$ is the cost function appearing in the variational formulation of multi-slice projector~\eqref{eq:MSdec}. \vspace{0.2cm}

Clearly $\projector[\Sps_\dimSps](\ztrue)\in\Sc$ because
\begin{align}
f(\projector[\Sps_\dimSps](\ztrue))
&= \sum_{j=1}^\dimSo \kparen{\lop(\SoPGbasis_j)- \aop(\projector[\Sps_\dimSps](\ztrue),\SoPGbasis_j)}^2\nonumber\\
&= \sum_{j=1}^\dimSo \kparen{\scap[\rieszr_j]{\ztrue}- \scap[\rieszr_j]{\projector[\Sps_\dimSps](\ztrue)}}^2 \nonumber\\
&= \sum_{j=1}^\dimSo \kparen{\scap[\rieszr_j]{\projector[\Sps_\dimSps][\perp](\ztrue)}}^2 \nonumber\\
& \leq \gamma^2 \kvvbar{\projector[\Sps_\dimSps][\perp](\ztrue)}^2\nonumber\\
&\leq \gamma^2 \Spwidth_\dimSps^2. \label{eq:defBSf}
\end{align}

Moreover, $\zrommult\in\Sc$.~This can be seen from the following arguments.~First, $\projector[\Sps_\dimSps](\ztrue)$ is a feasible point for problem \eqref{eq:MSdec}, that is 
\begin{align}
\dist[{\projector[\Sps_{\dimSps}](\ztrue)}]{\Sps_\idxss}\leq \Spswidth_\idxss \mbox{ for $\idxss = 0\ldots \dimSps$.} 
\end{align}
Indeed, rewriting $\ztrue$ as
\begin{align}
\ztrue 
&= \sum_{j=1}^\dimSps \scap[\Spsbasis_j]{\ztrue} \Spsbasis_j + \zz, 
\end{align}
where $\zz\in\Spsp_\dimSps$, we have
\begin{align}
\Spswidth_\idxss
&\geq \dist[{\ztrue}]{\Sps_\idxss}\nonumber\\
&=  \kvvbar{{\projector[\Sps_{\idxss}][\perp](\ztrue)}}\nonumber\\
&= \kvvbar{\sum_{j=\idxss+1}^\dimSps \scap[\Spsbasis_j]{\ztrue} \Spsbasis_j + \zz}\nonumber\\
&= \sqrt{\kvvbar{\sum_{j=\idxss+1}^\dimSps \scap[\Spsbasis_j]{\ztrue} \Spsbasis_j}^2 + \kvvbar{\zz}^2}\nonumber\\
&\geq \kvvbar{\sum_{j=\idxss+1}^\dimSps \scap[\Spsbasis_j]{\ztrue} \Spsbasis_j}\nonumber\\
&=\kvvbar{{\projector[\Sps_{\idxss}][\perp]({\projector[\Sps_{\dimSps}](\ztrue)})}}\nonumber\\
&=\dist[{\projector[\Sps_{\dimSps}](\ztrue)}]{\Sps_\idxss}.\label{eq:Pztruefeas}
\end{align}
The first inequality follows from our initial assumption $\ztrue\in \cap_{\idxss=0}^\dimSps \slice_{\idxss}$. The third equality is true because $\zz\in\Spsp_\dimSps$. 
 Now, since $\zrommult$ is a minimizer of $f(\vech)$ over the set of feasible points, we have $f(\zrommult)\leq f(\projector[\Sps_\dimSps](\ztrue))\leq \gamma^2 \Spwidth_n^2$ and therefore $\zrommult\in\Sc$. \vspace{0.2cm}

We finally show that $\zrommult\in\Sc$ and $\projector[\Sps_\dimSps](\ztrue)\in\Sc$ implies $\projector[\Sps_\dimSps](\ztrue)-\zrommult\in\feasSd_1$. 
 Let us first note that, if $\vech\in\Sps_\dimSps$, the cost function $f(\vech)$ can be rewritten as: 
 \begin{align}
f(\vech) 
&= \sum_{j=1}^\dimSo \kparen{\lop(\SoPGbasis_j) - \aop(\vech,\SoPGbasis_j)}^2\nonumber\\
&= \sum_{j=1}^\dimSo \kparen{\scap[\rieszr_j]{\ztrue} - \scap[\rieszr_j]{\vech}}^2,\nonumber\\
&= \sum_{j=1}^\dimSo \kparen{\scap[\rieszr_j^*]{\ztrue} - \scap[\rieszr_j^*]{\vech}}^2,\nonumber\\
&= \sum_{j=1}^\dimSps \kparen{\scap[\rieszr_j^*]{\ztrue} - \svdecr_j \scap[\Spsbasis_j^*]{\vech}}^2 + \sum_{j=\dimSps+1}^\dimSo \scap[\rieszr_j^*]{\ztrue}^2,\label{eq:exprfun}
\end{align}
where the third equality follows from the fact that $\kbrace{\rieszr_j}_{j=1}^\dimSo$ and $\kbrace{\rieszr_j^*}_{j=1}^\dimSo$ differ up to an orthonormal transformation; the last equality is a consequence of \eqref{eq:orthproperty} and the fact that $\vech\in\Sps_\dimSps$ by hypothesis. \vspace{0.2cm}

 We note that  $\projector[\Sps_\dimSps](\ztrue)-\zrommult$ can be written as $\sum_{j=1}^{\dimSps} \beta_j \Spsbasis_j^*$ by setting $\beta_j = \scap[\Spsbasis_j^*]{\projector[\Sps_\dimSps](\ztrue)} - \scap[\Spsbasis_j^*]{\zrommult}$.
 Therefore, we have
 \begin{align}
\sum_{j=1}^{\dimSps} \svdecr_j^2 \beta_j^2
=& \sum_{j=1}^{\dimSps} \kparen{\svdecr_j \scap[\Spsbasis_j^*]{\projector[\Sps_\dimSps](\ztrue)} - \svdecr_j\scap[\Spsbasis_j^*]{\zrommult}}^2, \nonumber\\
=& \sum_{j=1}^{\dimSps} \Bigl( \svdecr_j\scap[\Spsbasis_j^*]{\projector[\Sps_\dimSps](\ztrue)}-\scap[\rieszr_j^*]{\ztrue} - \svdecr_j \scap[\Spsbasis_j^*]{\zrommult}+\scap[\rieszr_j^*]{\ztrue}\Bigr)^2, \nonumber\\
\leq&  2\sum_{j=1}^{\dimSps} \kparen{\svdecr_j \scap[\Spsbasis_j^*]{\projector[\Sps_\dimSps](\ztrue)}-\scap[\rieszr_j^*]{\ztrue}}^2 +2 \sum_{j=1}^{\dimSps}   \kparen{\svdecr_j \scap[\Spsbasis_j^*]{\zrommult}-\scap[\rieszr_j^*]{\ztrue}}^2,\nonumber\\
\leq&  2 f(\projector[\Sps_\dimSps](\ztrue)) + 2 f(\zrommult),\nonumber\\
\leq& 4 \gamma^2 \Spwidth_\dimSps^2,\nonumber
\end{align}
where the first inequality follows from the standard inequality $(a+b)^2\leq2 (a^2+b^2)$, the second from \eqref{eq:exprfun},  and the last one from the fact that $\zrommult\in\Sc$ and $\projector[\Sps_\dimSps](\ztrue)\in\Sc$.

\subsubsection{Definition of $\feasSd_2$}\label{sec:defD2}
Let 
\begin{align}
\delta_j &= \tmpvar_j + \hat{\tmpvar}_j,
\end{align}
where
\begin{align}
\tmpvar_j          &= \sum_{i=1}^\dimSps \kvbar{x_{ij}} \Spwidth_{i-1},\nonumber\\
\hat{\tmpvar}_j &= \sum_{i=1}^\dimSps \kvbar{x_{ij}} \Spswidth_{i-1},
\end{align}
and the $x_{ij}$'s are the elements of the matrix $\X$ appearing in the SVD decomposition \eqref{eq:svd_Gr}. We define $\feasSd_2$ as
\begin{align}
\feasSd_2 = \kbrace{\feasd = \sum_{j=1}^{\dimSps} \beta_j \Spsbasis_j^* : \kvbar{\beta_j}\leq  \tmpvar_j}.
\end{align}
We show hereafter that $\projector[\Sps_\dimSps](\ztrue)-\zrommult\in \feasSd_2$. \vspace{0.2cm}

We first note that if $\vech$ is feasible for problem \eqref{eq:MSdec}, we must have
\begin{align}
\kvbar{\scap[\Spsbasis_j^*]{\vech}}\leq \hat{\tmpvar}_j.\label{eq:boundscapvj1}
\end{align}
Indeed, if $\vech$ is feasible, the constraint $\dist[{\vech}]{\Sps_\idxss}\leq \Spswidth_\idxss$ simply writes as
\begin{align}
\sum_{j=\idxss+1}^\dimSps \scap[\Spsbasis_j]{\vech}^2 \leq \Spswidth_\idxss^2.\nonumber
\end{align}
In particular, this implies that 
\begin{align}
\kvbar{\scap[\Spsbasis_{\idxss+1}]{\vech}} \leq \Spswidth_\idxss. \nonumber
\end{align}
Using the fact that
\begin{align}
\Spsbasis_j^* &= \sum_{\idxss=1}^\dimSps x_{\idxss j} \Spsbasis_\idxss,\nonumber
\end{align}
we obtain \eqref{eq:boundscapvj1}. 
In a similar way, we can find that
\begin{align}
\kvbar{\scap[\Spsbasis_j^*]{\projector[\Sps_\dimSps](\ztrue)}}\leq \tmpvar_j,\label{eq:boundscapvj2}
\end{align}
by using the fact that $\dist[{\projector[\Sps_\dimSps](\ztrue)}]{\Sps_\idxss}\leq \Spwidth_\idxss$ from \eqref{eq:Pztruefeas}. \vspace{0.2cm}

Let us now show that $\projector[\Sps_\dimSps](\ztrue)-\zrommult\in \feasSd_2$. 
We first note that $\projector[\Sps_\dimSps](\ztrue)-\zrommult$ can be written as $\sum_{j=1}^{\dimSps} \beta_j \Spsbasis_j^*$ by setting $\beta_j = \scap[\Spsbasis_j^*]{\projector[\Sps_\dimSps](\ztrue)} - \scap[\Spsbasis_j^*]{\zrommult}$. This leads to
\begin{align}
\kvbar{\beta_j}
&=\kvbar{\scap[\Spsbasis_j^*]{\projector[\Sps_\dimSps](\ztrue)} - \scap[\Spsbasis_j^*]{\zrommult}},\nonumber\\
& \leq \kvbar{\scap[\Spsbasis_j^*]{\projector[\Sps_\dimSps](\ztrue)}} + \kvbar{\scap[\Spsbasis_j^*]{\zrommult}},\nonumber\\
& \leq \hat{\tmpvar}_j + \tmpvar_j = \delta_j, \nonumber
\end{align}
where the last inequality follows from \eqref{eq:boundscapvj1} and \eqref{eq:boundscapvj2}.

\subsection{Expression of $ \sup_{\feasd\in\feasSd}\kvvbar{\feasd}^2$}

 We consider the following problem: 
\begin{align}
 \sup_{\feasd\in\feasSd}\kvvbar{\feasd}^2&= \sup_{\boldsymbol{\beta}} \kvvbar{\boldsymbol{\beta}}^2\ \mbox{ subject to  }
\left\{
\begin{array}{l}
\sum_{j=1}^{\dimSps} \svdecr_j^2 \beta_j^2\leq 4 \gamma^2 \Spwidth_\dimSps^2\\
\kvbar{\beta_j} \leq \delta_j  
\end{array}
\right. .\label{eq:optimizationbound}
\end{align}
If $\sum_{j=1}^{\dimSps} \svdecr_j^2 \delta_j^2\leq 4 \gamma^2 \Spwidth_\dimSps^2$, the first constraint in \eqref{eq:optimizationbound} is always inactive and the solution simply reads
\begin{align}
 \sup_{\feasd\in\feasSd}\kvvbar{\feasd}^2 = \sum_{j=1}^{\dimSps} \delta_j^2. 
\end{align}
If $\sum_{j=1}^{\dimSps} \svdecr_j^2 \delta_j^2\geq 4 \gamma^2 \Spwidth_\dimSps^2$, 
the solution of \eqref{eq:optimizationbound} is given by 
\begin{align}
 \sup_{\feasd\in\feasSd}\kvvbar{\feasd}^2 &= \sum_{j=\ell+1}^\dimSps \delta_j^2+ \rho\, \delta_\ell^2,\label{eq:bornesup}
\end{align}
where $\ell$ is the largest integer such that
\begin{align}
\sum_{j=\ell}^{\dimSps} \svdecr_j^2 \delta_j^2\geq 4 \gamma^2 \Spwidth_\dimSps^2,
\end{align}
and $\rho\in[0,1]$ is defined as
\begin{align}
\rho \svdecr_\ell^2 \delta_\ell^2 +\sum_{j=\ell+1}^{\dimSps} \svdecr_j^2 \delta_j^2= 4 \gamma^2 \Spwidth_\dimSps^2. 
\end{align}
This can be seen by verifying the optimality condition of problem \eqref{eq:optimizationbound}. 
We note that problem \eqref{eq:optimizationbound} is the same (up to some constants) to the one considered in \cite[Section 3.1]{Binev2015Data}. The solution \eqref{eq:bornesup} is therefore similar, up to some different constants, to the one obtained in that paper. 

\section*{Acknowledgement}

The authors thank the ``Agence nationale de la recherche'' for its financial support through the Geronimo project (ANR-13-JS03-0002).

\bibliographystyle{ieeetr}
\bibliography{/Users/cherzet/Dropbox/Personnel/Biblio/MyBibli.bib}


\end{document}